\documentclass[12pt, oneside]{article}   	
\pdfoutput=1
\usepackage{geometry}                		
\usepackage{enumerate}
\usepackage{graphicx}				
\usepackage{subcaption}
\usepackage{amssymb}
\usepackage{amsmath}
\usepackage{amsthm}
\usepackage{amsfonts}
\usepackage{proof}
\usepackage{ifthen}
\usepackage{centernot}
\usepackage{tikz-cd}
\usepackage{lineno}
\usepackage{multirow}
\usepackage{booktabs}
\usepackage[authoryear]{natbib}

\newcommand{\be}{\begin{equation}}
\newcommand{\ee}{\end{equation}}

\everymath{\displaystyle}

\theoremstyle{definition}

\theoremstyle{remark}

\newtheorem
{notation}{Notation}
\newtheorem
{note}{Note}
\title{Critical speeding up as an early warning signal of regime switching.}
\author{Mathew Titus$^{1,2}$, Zach Gelbaum$^{1,2}$, James R. Watson$^{1}$}

\begin{document}
\maketitle

$^1$College of Earth, Ocean and Atmospheric Sciences, Oregon State University, Corvallis, Oregon, USA\\
$^2$The Prediction Lab, Corvallis, Oregon, USA\\



\clearpage

\begin{abstract}
The use of critical slowing down as an early warning indicator for regime switching in observations from stochastic environments and noisy dynamical models has been widely studied and implemented in recent years. Some systems, however, have been shown to avoid critical slowing down prior to a transition between equilibria, e.g. \citep{ditlevsen2010wishfulthinking}. Possible explanations include non-smooth potential driving the dynamic \citep{hastings2010nowarning} or large perturbations driving the system out of the initial basin of attraction. In this paper we discuss a phenomenon analogous to critical slowing down, where a slow parameter change leads to a high likelihood of a regime shift and creates signature warning signs in the statistics of the process's sample paths. In short, if a basin of attraction is compressed under a parameter change then the potential well steepens, leading to a drop in the time series' variance and autocorrelation; precisely the opposite warning signs exhibited by critical slowing down. This effect, which we call ``critical speeding up,'' is demonstrated using a simple ecological model exhibiting an Allee effect. The fact that both dropping and rising variance / autocorrelation can indicate imminent state change should underline the need for reliable modeling of any empirical system where one desires to forecast regime change.
\end{abstract}

When studying time series data for dynamical systems which exhibit critical transitions, that is sudden changes in equilibrium behavior, a widely used early warning signal for an oncoming transition is critical slowing down (CSD) \citep{scheffer2009early,ellis2012ewsmethods}. This signature for the system being at risk of a large transition is based on the theory of stochastic dynamical systems \citep{kuehn2013mathematical}, but has been observed in a variety of empirical tests, both in nature and in the laboratory (see \citep{carpenter2006risingvariance, vanbelzen2017marshes}, \citep{wen2018one} and references therein). At its core, CSD assumes that the stochastic process $X$ experiences a smooth potential $V_t$ which is varying slowly in time, $dX_t = -\nabla V_t(X_t)dt + \xi_t$ where $\xi$ is some random process.  If the potential nears a bifurcation point, such as a pitchfork or fold bifurcation, then the shape of $V$ around its equilibrium necessarily flattens out as its minimum (stable point) becomes a degenerate critical point.  The lessening or loss of local curvature, responsible for the mean-reverting property within the basin of attraction, means that excursions of $X$ away from its stable point grow both in extent and length of time. The concomitant increase in variance and autocorrelation of the sample path of $X$ is what we refer to as critical slowing down. We expand on this phenomenon in the next section. 

The catalogued examples where CSD is observed are often carefully controlled laboratory settings \citep{kramer1985bistable} or models of natural systems subject to relatively small perturbations \citep{carpenter2006risingvariance,dakos2014critical}, but many empirical observations of dynamical systems fail to exhibit critical slowing down prior to making a change of regime \citep{ditlevsen2010wishfulthinking,hastings2010regime,rozek2017evidence}. 
That critical state transitions appear to commonly occur in nature away from a bifurcation in the underlying governing dynamic speaks to the value in devising early warning indicators for other `high-risk' situations where a regime shift may be becoming increasingly likely.
This lack of CSD prior to transition has been discussed in the literature recently, with explanations including large exogenous perturbations\citep{boettiger2013largedeviations,ditlevsen2010wishfulthinking} or the potential $V$ lacking smoothness \citep{hastings2010nowarning}; these translate, respectively, to a sudden change in the governing equation (a fundamental change in the statistics of $\xi$) or a departure from the classical setting of modeling the dynamics with a system of partial differential equations and perturbing the smooth dynamical model with the noise source $\xi$.

In the present paper we study an alternative culprit which retains both smoothness of the potential and the (quasi-)stationarity of the system. We assume the usual setting for CSD: a dynamical system governed by a slowly-changing potential $V_t$, and a fixed noise term, $\xi_t = a B_t$ for some $a > 0$ and $B$ a brownian motion. If $V$ experiences a bifurcation the degeneracy of $\nabla V_t$ will create a slowing effect, increasing both variance and autocorrelation. 
However, we consider instead potentials where the long-term evolution of $V$ narrows the basin of attraction that the process is initiated within, avoiding bifurcations. 
Given a potential well with fixed height but shrinking width leads to a `rattling' effect, wherein the stabilizing gradient grows, punishing excursions away from the stable point (minimum of $V$) and shortening their extent. 
This manifests itself in decreased variance and decreased autocorrelation of the sample path of $X$ or as we call it ``rattle''.
These signs indicate that a stochastic transition is becoming more likely, giving us an early warning of a stochastic transition, rather than a transition compelled by bifurcation. (Note that this is not at odds with the conclusions of \citep{boettiger2013largedeviations}, for the transition is not purely due to noise.)

We further show that if the narrowing of the potential is uniform (a linear rescaling of the spatial axis) then the effect on the exit time distribution is equivalent to rescaling the time variable. 
That is, in terms of the risk of escaping the basin of attraction within the next $T$ units of time, narrowing the potential is the same as speeding up the process's evolution, and from this standpoint it is clear that such a transformation increases the chance of regime switching. 
This is the motivation for our term ``critical speeding up'' (CSU), as the time until a critical transition is shortened.

In the next section we give a simple mathematical description of both critical slowing down and critical speeding up. As an example of CSU occurring in a simple ecological model, Section \ref{s:pop_model} introduces a population-level model for a species with an Allee effect (i.e. when individual fitness and/or the population growth rate is a function of population density). We assume that for this species that the population growth rate is proportional to the population density; then under a shrinking habitat the system begins to rattle (see Figure \ref{fig:csu}) and exhibit CSU effects as risk of extinction grows. In the final section we discuss the broader implications of this effect, in particular that the results presented here stress the need to implement a faithful model before relying on either CSD or CSU as an early warning signal.

\section{Mathematical background: CSD and CSU}
\label{s:math}

Suppose for simplicity that we are dealing with a one-dimensional process $X = \{ X_t : t \geq 0 \}$ where $X_t$ may be a population size, or a concentration of a chemical, or a fraction of people who subscribe to a given belief, at time $t$. The systems we are interested in allow $X$ to be modeled by a stochastic differential equation (SDE) as it generally obeys some dynamical law, but is subject to exogenous perturbations or uncertainty in measurement:
\begin{equation}
\label{sde}
dX_t = b(X_t,t)dt + a(X_t,t) dB_t.
\end{equation}
Here $b \in C^{\infty}$ is the (smooth) drift function, describing the deterministic differential equation $X$ would obey if $a \equiv 0$; the noise is modeled by the product of the coefficient $a \in C^{\infty}$ and $B = \{ B_t : t \geq 0 \}$, a one-dimensional brownian motion. 
Since we are working in one-dimension, a potential function exists
\begin{equation}
V(x,t) = -\int_0^x b(x,t) dx
\end{equation}
so that (\ref{sde}) becomes $dX_t = -\nabla V(X_t,t) dt + a(X_t,t) dB_t$.

If the process is initiated in a basin of attraction which, over time, is bounded on the left (right) by $x_l(t)$ (respectively, $x_r(t)$), i.e. $X_0 = x \in \left(x_l(0), x_r(0) \right)$, we write $T(x)$ for the exit time of the process from that basin:
\begin{equation}
\label{exit_time}
T(x_0) = \inf_{t \geq 0} \{ X_t \not\in (x_l(t), x_r(t)) \vert X_0 = x \}.
\end{equation}
In the sequel we will use the potential
\begin{equation}
\label{toy_model}
V(x) = \beta^3 x^3 - \alpha \beta x
\end{equation}
as a testbed for both phenomena.

\subsection{Equivalence of narrowing and time change}
\label{ss:time_change}

Suppose the process $X = \{ X_t \}$ obeys the SDE $dX_t = -\nabla V(X_t)dt + a dB_t$, then if the spatial scale is compressed by a factor $k > 1$ so that the potential defining the dynamic becomes $\hat{V}(x) := V(kx)$, then (leaving the diffusion term unchanged) we have a second process, $\hat{X} = \{ \hat{X}_t \}$, which evolves under the influence of the narrowed potential $\hat{V}$, and so experiences at point $x$ a drift $-\nabla \hat{V}(x) = -k \nabla V (kx)$:
\begin{equation}
\label{x_hat}
d\hat{X}_t = - k \nabla V(k \hat{X}_t) + a dB_t. 
\end{equation}

Let us define $Y_t := k \hat{X}_{t/k^2}$. Then from (\ref{x_hat}) and the identity $\sqrt{s} B_t \sim B_{st}$ we have
\begin{eqnarray*}
k d\hat{X}_t & = & - k^2 \nabla V(k \hat{X}_t) dt + a k dB_t, \\
d Y_{k^2 t} & = & - \nabla V(Y_{k^2 t}) d(k^2 t) + a dB_{k^2 t},
\end{eqnarray*}
and so taking $s := k^2 t$ we find $Y = \{ Y_s : s \geq 0 \}$ obeys $dY_s = -\nabla V(Y_s) ds + adB_s$, the same equation as $X$.

Since $Y = \{ Y_s \}$ obeys the same dynamic as $X = \{ X_t \}$, the probability of $Y$ exiting the basin $(x_l, x_r)$ by some time $s = k^2 t < k^2 \tau$ (when initiated from $x$) is equal to $p_{k^2 \tau}(x)$. Define $\hat{T}(x) = \inf_{t \geq 0} \left\{ \hat{X}_t \not\in (x_l/k, x_r/k) \right\} = \inf_{s \geq 0} \left\{ Y_s \not\in (x_l, x_r) \right\}$ and $\hat{p}_{\tau}(x) = P\left( \hat{T}(x) < \tau \right)$. 
Then finally we have the following relation between exit time distributions:
\begin{equation}
\label{time_change}
 \hat{p}_{\tau}(x) = p_{k^2 \tau}(x)
\end{equation}
Accordingly, contracting space by a factor $k > 1$ leads to the same exit time statistics (i.e. critical transition rates) as contracting time by a factor of $k^2$. Effectively, a uniform narrowing of $V$ can be interpreted as speeding up the evolution of $X$, increasing the rate of extreme events.

We note that the above argument also applies for $0 < k < 1$, in which case the potential widens, and has the same statistics as a copy of the process which evolves a factor of $k^2$ slower; this slowing effect is the origin of critical slowing down.

\subsection{CSD and CSU}

In this subsection we discuss the fundamentals of CSD; it is not intended to be an exhaustive treatment and readers should turn to the literature for full details and generality \citep{kuehn2013mathematical}. As above, we assume that the process of interest is one-dimensional and described by the SDE $dX_t = - \nabla V_t(X_t) dt + a dB_t$, where $V$ is smooth, $a$ is a positive constant, and $B$ is a brownian motion. As we are not interested in deriving the explicit estimates of CSD using normal forms and slow-fast systems theory \citep{berglund2006noise, kuehn2013mathematical}, but rather the qualititative features, we approximate the potentials with the lowest order term in their Taylor expansions. This simplifies the dynamics to those of the Ornstein-Uhlenbeck processes or brownian motions with drift.
We also assume that the process originates within a basin of attraction, $X_0 \in(x_l(0), x_r(0)) \subset \mathbb{R}$, and that the (smooth) functions $x_l(t), x_r(t)$ are defined so that the basin of attraction is given by $(x_l(t), x_r(t))$ for all future times that it exists.

Near the unique minimum of $V$ within the basin $(x_l(t),x_r(t))$ the second-order approximation of $V$ can be used in place of the true potential to define a process whose dynamics are close to those of $X$ while it remains near equilibrium. Since we are interested in the system (\ref{toy_model}), we consider the modified equation
\begin{equation}
\label{approx_sde}
dX_t = - 2 \beta^2 \sqrt{ 3 \alpha } X_t dt + a dB_t
\end{equation}
which simply follows from the Taylor expansion of $V$ about the point $x_0 = \sqrt{\alpha / 3 \beta^2}$.

As mentioned above, this process is an Ornstein-Uhlenbeck (OU) process so its characteristics are well-known.
In particular, the $s$-autocorrelation function defined by 
$$\textrm{Cov}(X_t,X_{t+s}) = \mathbb{E}\left[ (X_{t+s} - \mathbb{E}X_{t+s})(X_t - \mathbb{E}X_t) \right] $$
is given for the Ornstein-Uhlenbeck process $dX_t = -bX_t dt + a dB_t$ by the formula
\begin{equation*}
\textrm{Cov}(X_t,X_{t+s}) = {a^2 \over 2 b} \left( e^{-sb} - e^{-(2t+s)b} \right).
\end{equation*}
From this it is immediate that the variance is given by 
$$\textrm{Var}(X_t) := \mathbb{E}\left[ \left( X_t - \mathbb{E}[ X_t ] \right)^2 \right] = {a^2 \over 2b} \left( 1 - e^{-2 b t} \right), $$ 
and so we have for (\ref{approx_sde}) that the autocorrelation with $s = 1$ is given by
\begin{equation}
\label{cov}
\textrm{Cov}(X_t,X_{t+1}) = {a^2 \over 4\beta^2 \sqrt{3\alpha}} e^{-2 \beta^2 \sqrt{3\alpha}} \left( 1 - \exp\left( -4\beta^2\sqrt{3\alpha} t \right) \right).
\end{equation}

Now, if we allow $\alpha$ to change in time so that it slowly approaches zero from above, a bifurcation occurs at $\alpha = 0$. Let us define $\alpha(t) := t^{-2}$; then $\alpha$ will approach 0 slowly enough that CSD is observable, and we see that the variance grows linearly in time, $\textrm{Cov}(X_t,X_{t+!}) = mt$ for $m = {a^2 (1 - \exp(-4\beta^2\sqrt{3}) ) \over 4\beta^2\sqrt{3} }$, and similarly for the autocorrelation.


On the other hand, if we do not approach the bifurcation point, fixing $\alpha > 0$ in time, and instead take $\beta$ to be increasing in time, say $\beta(t) = t$ and restrict to $t \geq 1$, then, again translating the potential's minimum to the origin, we have 
$$ \nabla V_t(x)  = 2 \sqrt{3\alpha} t^2 x = \nabla V_0(tx), $$
demonstrating the speeding up incurred by narrowing the potential by a factor $t$. 

We may again apply (\ref{cov}) to recover the variance and autocorrelation. One sees that for this process
$$ \textrm{Var}(X_t) = {a^2 \over 4\sqrt{3 \alpha} t^2} \left( 1 - e^{-4\sqrt{3\alpha} t^3} \right) = O(t^{-2}), $$
and likewise the 1-step autocorrelation decays like $t^{-2}$.
While we have restricted our attention to a rather specific example, this decay of autocorrelation and variance, which we call critical speeding up, is characteristic of contracting potentials.

\section{Observation of CSU in a population model}
\label{s:pop_model}

Here we introduce a population model that exhibits CSU before a critical state transition (population collapse). We assume a strong Allee effect, proportional to the territory size; this could model a species which has displaced its competitors, benefitting from a cooperative advantage against another species. We also assume that the reproduction rate is a function of population density, rather than total population. These features create an attracting region for the fixed point at zero, so that after passing below a critical threshold the species will become extinct (discounting a restoring perturbation).

This model is described by the following equation
\begin{equation}
\label{allee_model}
{dx \over dt} = {r \over \beta} x \left( {x \over \beta A} - 1 \right) \left( 1 - {x \over \beta C} \right), ~~ \beta \in [0,1],
\end{equation}
where all parameters and the population size $x$ are assumed to be nonnegative. 

We interpret $C$ as the maximum carrying capacity of the species, say if they control $100\%$ of their potential territory, while $\beta C$ is the carrying capacity when a fraction $\beta$ is controlled. The product $\beta A$ represents the strength of the Allee effect; suppose for example that the population only controls one quarter of the contested territory ($\beta = 0.25$) rather than half ($\beta = 0.5$), then the Allee effect is weakened as their competitors have control of a greater swath of the environment, reducing the degree to which the individuals are able to benefit from cooperation. Finally, $r$ is chosen so that a population of size $x$ reproduces at a rate $rx$ when distributed over the entire territory. When restricted to a fraction $\beta$ of the total land, the population density increases by a factor of $1/\beta$. Hence, the reproductive rate is given by ${r \over \beta} x$.

As a final addition, we assume the model originates at its positive equilibrium ($X_0 = \beta C$) and add the noise term $\alpha B$ to model exogenous perturbations (change in population through fatal accidents, immigration from another distant population, twins, etc.), and write $X = \{ X_t : t \geq 0 \}$ for the population process:
\begin{equation}
\label{allee_sde}
dX_t = {r \over \beta} X_t \left( {X_t \over \beta A} - 1 \right) \left( 1 - {X_t \over \beta C} \right) dt + \alpha dB_t. 
\end{equation}

\begin{figure}
\centering
\includegraphics[width = 5.55in, height = 5.0in]{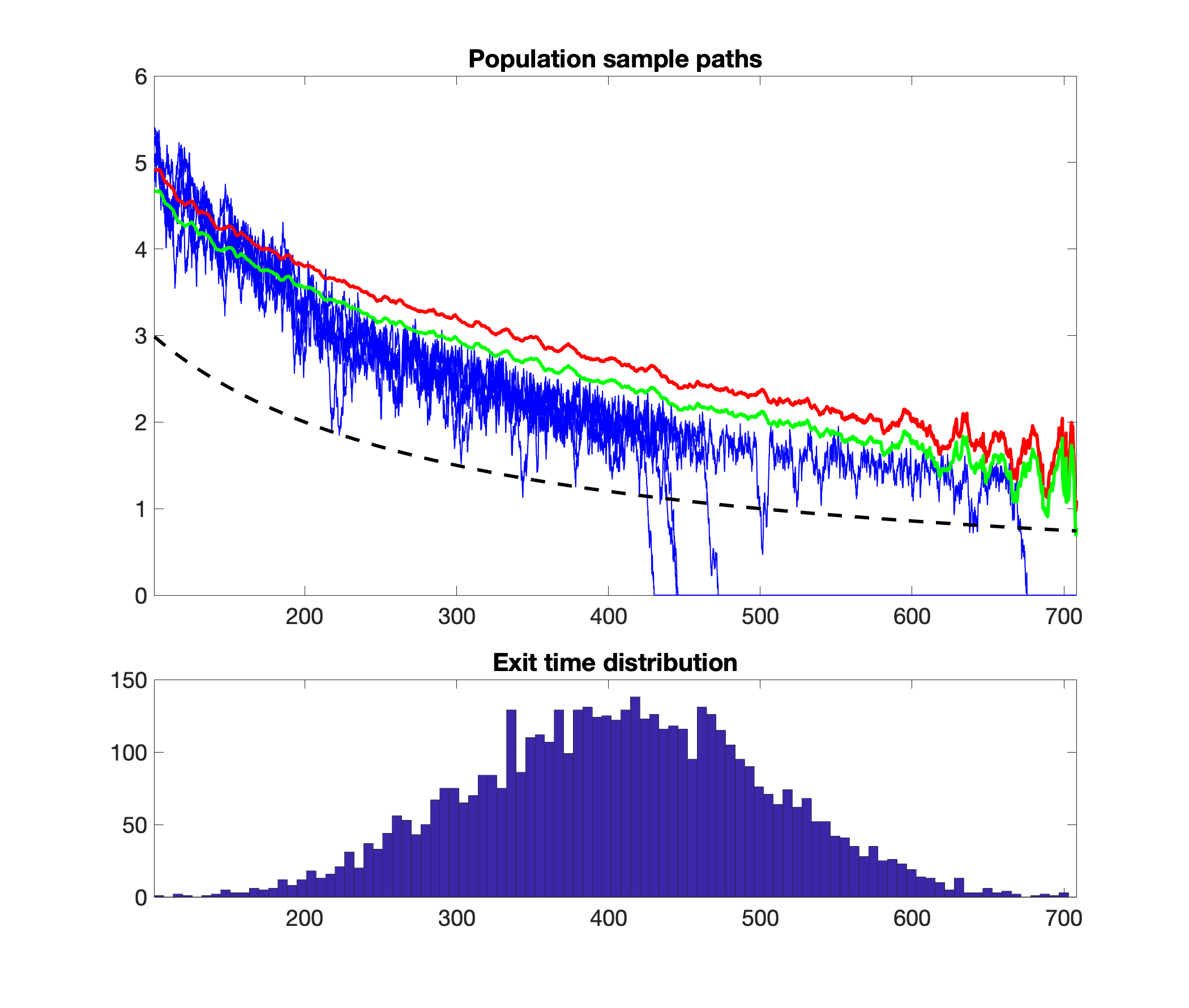}
\caption{{\footnotesize The top panel displays in blue five of the 5000 calculated sample paths of $X_t$ as $\beta$ decreases from 4 to roughly $1/2$. The mean variance of the surviving sample paths, calculated over a rolling window of length 10 is shown in red; the mean 1-step autocorrelation of the surviving paths within the same window is plotted in green. NB: On the far right of the figure the low number of surviving samples leads to a greater variance in these averages, as the law of large numbers breaks down. The dashed black line denotes $A\beta(t)$, the unstable equilibrium point dividing the extinction regime from the survival regime. \\
In the lower panel we display a histogram of the time required for $X_t$ to enter the basin of attraction of $x = 0$ (extinction regime).}}
\label{fig:csu}
\end{figure}

We study this model's vulnerability to extinction events under the increasing stress of an encroaching competing species, which decreases the available territory, driving $\beta$ towards zero:
\begin{equation}
\beta(t) = {4 \over 1 + 0.01 t}.
\end{equation}
In Figure \ref{fig:csu} we see the effect of slowly decreasing $\beta$ over 5000 independent trials; the lower panel shows the distribution of the time of collapse, which we define as the first time the system's population $X_t$ drops below $\beta(t) A$, as this separates the two dynamical regimes (extinction and survival). For all of these trials we fix $\alpha = 0.22$, $C = 2.5$, $A = 1.5$ and $r = 1$. Also pictured in Figure \ref{fig:csu} are the variance and autocorrelation functions (calculated over a rolling window of length 10), after averaging over all sample paths which have yet to exit the basin of attraction. The hallmarks of CSU are evident in the time series; as $\beta$ drops, the equilibrium points ($C\beta(t)$, $A\beta(t)$, and $0$) gather toward to origin, and both the variance and autocorrelation of the series are reduced. Despite the variance and autocorrelation dropping by over 50\%, none of the 5000 trials survived past 871 units of time (surviving sample path not pictured).

\begin{figure}
\centering
\includegraphics[scale = 0.25]{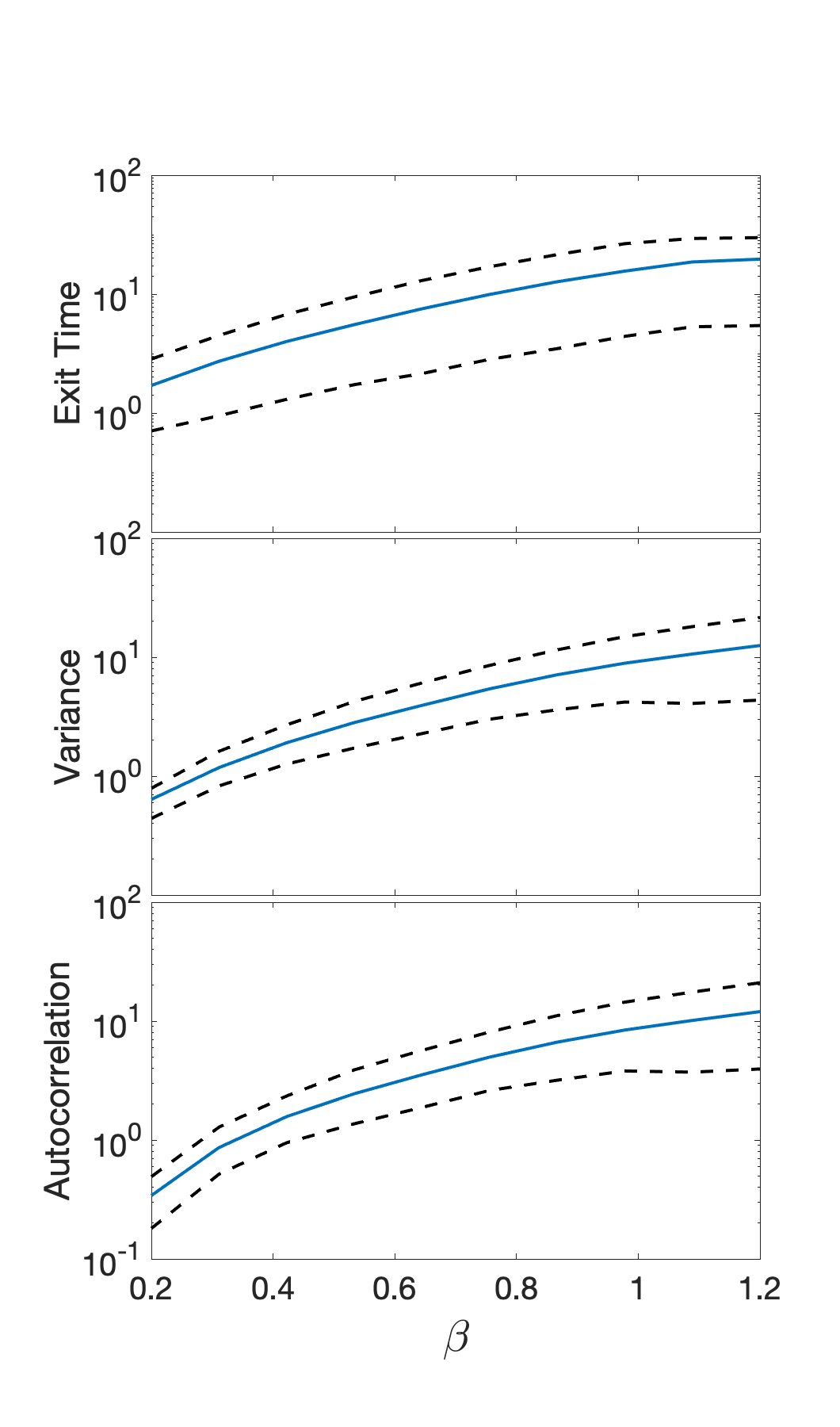}
\caption{{\footnotesize Here we give semi-log plots of the mean, 5th percentile, and 95th percentile for three statistics defining Critical Speeding Up. Top: exit time of $X_t$ from the basin of attraction $(A\beta, \infty)$. Middle: variance of $X_t$ over the 10 time units preceding its exit from the basin. Bottom: 1-step autocorrelation of $X_t$ over the 10 time units preceding its exit from the basin.}}
\label{fig:beta}
\end{figure}

As another verification of the theory in subsection \ref{ss:time_change}, we also run $2000$ repeated simulations of $X_t$ starting from $X_0 = \beta C$ for $\beta$ equal to a fixed value between 0.2 and 1.2 (ten equally distributed values were used, making for 20,000 independent trials). Statistics of the results are plotted in Figure \ref{fig:beta}. It is shown that the time for the population $X_t$ to exit $(\beta A, \infty)$, the basin of attraction about $\beta C$, decreases with decreasing $\beta$ as expected from the previous analysis, i.e. the risk of a critical transition grows as the potential narrows. However, the calculation of the covariance and autocorrelation (see Figure \ref{fig:beta}) show that the system becomes more brittle as $\beta$ drops towards $0.2$, and critical speeding up gives the precursor signal for population collapse. Note that all values are plotted on a semi-log scale.

\section{Discussion}

Using stochastic dynamical systems theory we have shown that noisy dynamical systems such as seen in nature, can not only fail to exhibit critical slowing down prior to a regime shift, they may actually exhibit a {\em speeding up} of their dynamic, with a decrease in both variance and autocorrelation. To a practitioner observing a system with a narrowing potential, a naive application of the theory of CSD would suggest the system is stabilizing as it begins to rattle, despite the impending critical transition. We used an example model of population growth to demonstrate this in practice, showing that a population experiencing territory loss and an Allee effect exhibits critical speeding up prior to collapse. We hope that this study underscores the necessity of having high-fidelity descriptions of the system dynamics before applying an early warning indicator such as CSD or CSU, a sentiment expressed also in \citep{boettiger2012quantifying}.

One might wonder just how pathological a system must be for it to compress a potential well $V$ under a slow parameter shift without significantly decreasing the depth of the well (leading to a bifurcation and CSD), nor increasing the depth of the well (restoring stability, despite the onset of CSU). We acknowledge that if the well deepens, the signal of CSU will be a false positive. However, this is analogous to a potential well broadening while deepening which leads to false positives under the CSD paradigm. Both signals offer an early warning under the appropriate conditions, but neither is a panacea for predicting regime change. Only with prior knowledge of the space of all possible evolutions of the underlying potential, $\{ V_t : t \geq 0 \}$, such as in the toy model (\ref{toy_model}), can we begin to map sample path statistics to high or low risk states with confidence.

As a second example of time series statistics belying the approach of a transition, one can imagine the model system (\ref{toy_model}) from Section \ref{s:math} evolving such that both $\alpha \to 0^+$ and $\beta \to \infty$. This would create both CSD and CSU effects in the time series. As these signals interfere with one another, simply observing the autocorrelation and variance of the time series may fail to find either a slowing down or a speeding up of the time series. This is another example of abrupt transitions occurring without a precursor slowing down, and may be responsible for some of the examples in the literature of critical transitions without CSD such as the regime switching in \citep{ditlevsen2010wishfulthinking}. At the least, this should call to question the idea that using CSD as an early warning indicator with no assumptions or prior knowledge of the system studied, a practice promoted elsewhere in the literature.

There are several reasons that we expect to find systems in nature that exhibit CSU. First, the model above describes a population governed by three simple rules, so it would be no surprise to find such a population within the earth's ecology. Further, suppose a sequence of failures must occur for a complex system (e.g. a power grid, social network or a full food-web) to experience catastrophe. If its subsystems become correlated, then the probability of a submodule failing (say, a generator shutting down, individuals disbanding, multiple species going extinct etc.) may be lowered due to the individuals' pooled resilience, but the number of failures required is lowered as well. In this case, the probability of an excursion of the state from stability ($X_0 = x_0$) to instability ($X_{T(x_0)}$) may be the same, but the number of events necessary (i.e. the distance in state space) has been reduced. This sort of thresholding between stability and instability is not uncommon in man-made systems, see for example \citep{reason1995error}.

As a third example, many engineered systems obey stochastic differential-algebraic equations (SDAEs), wherein the noisy dynamical system is required to satisfy some algebraic constraints \citep{schein1998sdaes}, so that for some function $f$ we have
\begin{eqnarray*}
dX_t &=& b(X_t, t) dt + a(X_t, t) dB_t, \\
f(X_t,t) &=& 0.
\end{eqnarray*}
Forcing the solution $X_t$ to remain near the zero set of $f(x,t)$ gives an overdetermined solution set; when solutions exist they are often ``stiff,'' meaning that the solution has a high-frequency component that can lead to sudden large-scale change in $X$. As the solution is subject to the condition $dX_t \cdot \nabla f(X_t,t) = 0$, we may reinterpret the system as having a strong restoring potential $V_f(x,t)$ such that solutions to $f(x,t) = 0$ minimize $V_f$, so that $X$ obeys an equation of the form $ dX_t = -\nabla(V + V_f) dt + dB_t$. Of course, the strength of $\nabla V_f$ suggests that $V_f$ is naturally quite narrow about its minima.

When studying a complex system with a variety of possible stable states, the notion of CSU may be useful in determining the likelihood of the system occupying a given equilibrium. Despite stabilizing drift that may be quite strong, if the basin of attraction about the stable point is prone to narrowing the system may not inhabit the basin for long periods of time. This principle may be useful in studying the time-dynamics and evolution of protein-folding, food-webs, or social structures \citep{zwanzig1997proteinfolding,ma2017foodchain,nekovee2007rumour}.

Finally, we mention that the principle of critical speeding up - that spatial compression of a basin of attraction can lead to regime switching as well as diminished autocorrelation and variance - holds in higher dimensions. One can also allow for systems which have stable limit cycles rather than equilibrium points. We do, however, relegate the proof of this fact to the larger literature on SDEs and future work.

\subsection*{Acknowledgments}
The author would like to acknowledge support from the DARPA YFA project N66001-17-1-4038 and George Hagstrom for suggesting the ecological model used here.
 

\bibliographystyle{achemso}
\bibliography{refs}

\end{document}